\documentclass[12pt,reqno]{amsart}
\usepackage{amsmath,amsfonts,  amssymb,amsthm,amscd}
\makeindex

\usepackage{amssymb}
\usepackage{graphics}
\usepackage[latin1]{inputenc}
\usepackage{epsf}

\usepackage{enumerate}
\usepackage{indentfirst}
\usepackage{euscript}
\usepackage{graphicx}
\usepackage{amsmath}
\usepackage{amsfonts}
\usepackage{amssymb}
\makeindex

 \theoremstyle{plain}
\newtheorem{proposition}{Proposition}
\newtheorem{lemma}{Lemma}

\newtheorem{remark}{Remark}
\theoremstyle{definition}

\newtheorem{theorem}{Theorem}

\newtheorem{example}[theorem]{Example}

\title[Principal Mean Curvature Foliations on Surfaces  ]{Principal Mean Curvature Foliations on Surfaces immersed in
${\mathbb R} ^4$}

\author{R. GARCIA, \;\;\;\; L.~F. MELLO \;\;\;
and  \;\;\;J. SOTOMAYOR }

 \thanks{ {T}his work is supported {CNP}q
\uppercase{G}rant 476886/2001-5 and  {PRONEX/FINEP/MCT - C}onv.
76.97.1080.00}

\begin{document}
\maketitle
\begin{abstract} {Here are studied qualitative properties of the families
of curves --foliations-- on a surface
immersed in  ${\mathbb R}^4$, along which it bends extremally
in the direction of the  mean  normal curvature vector.
 Typical singularities and cycles are described, which provide
 sufficient conditions, likely to be also necessary,   for the
 structural stability of the configuration of such foliations
 and their singularities, under small $C^3$ perturbations of the immersion.
   The conditions are expressed in terms of Darbouxian type of  the normal
   and umbilic singularities, the hyperbolicity of cycles,  and the asymptotic
    behavior of singularity separatrices and other typical curves of the foliations.
     They extend those given by Gutierrez and Sotomayor in 1982 for  principal foliations
      and umbilic points  of  surfaces  immersed in ${\mathbb R}^3$. Expressions for the
      Darbouxian conditions and for the hyperbolicity, calculable in terms of the
      derivatives of the immersion at singularities and cycles, are provided.
      The connection of the present  extension from ${\mathbb R}^3$ to ${\mathbb R}^4$
        to other pertinent ones as well as some problems left open in this
        paper are proposed at the end. }
\end{abstract}
\section{Position of the Problem, Results and Examples}\label {sec:1}

Let ${\mathbb M}^2$ be a 2-dimensional, compact, oriented, smooth manifold. Denote by ${\mathcal I}^{r,s}$ the space of immersions $\alpha$ of class $C^r$  of ${\mathbb M}^2$  into ${\mathbb R} ^4$,  with the $C^s$ topology, $r\geq s$.
${\mathbb M}^2$  and  ${\mathbb R} ^4$ are endowed with a once for all fixed orientation.

The projections of the pullback, ${\alpha}{^*} ({\mathbb R} ^4)$,    of the tangent bundle of  ${\mathbb R} ^4$ onto the tangent, ${\mathbb T}{\mathbb M}^2$, and  normal, ${\mathbb N}_{\alpha}$,  bundles of an immersion $\alpha$ will be  denoted respectively by  $\Pi_{\alpha,t}$ and $\Pi_{\alpha,n}$. These vector bundles are endowed with the standard metrics induced by the Euclidean one, $<,>$, in ${\mathbb R} ^4$.

Denote by $H_{\alpha}$ the {\it normal mean curvature vector field} of ${\alpha}$, which is defined by   centers of the field of {\it ellipses of normal curvature} of ${\alpha}$. For  any  positive orthonormal tangent frame $\{e_1, e_2\}$, it holds that $H_{\alpha}= \Pi_{\alpha,n}(De_1(e_1)+De_2(e_2))/2$. See  Little  \cite{li} and Garcia and Sotomayor  \cite{axial}.

The {\it normal singularities} of $\alpha$, defined by the  zeros of  $H_{\alpha}$,  will be denoted by ${\mathcal S}_{n,\alpha}$.
For  generic immersions,  when non-empty, ${\mathcal S}_{n,\alpha}$ is a discrete set of points at which $H_{\alpha}$ is  transversal to the zero section of ${\mathbb N}_{\alpha}$.    See Little   \cite{li} and Mello   \cite{lf}.

The {\it unit   normal mean curvature  vector field} of $\alpha$,  $N_\alpha =\frac{H_{\alpha}}{|H_{\alpha}|}$, is defined on the complement of  ${\mathcal S}_{n,\alpha}$.
The unit  vector field $B_\alpha = e_1 \wedge e_2 \wedge N_\alpha $ is well defined for  any  positive orthonormal tangent frame $\{e_1, e_2\}$ and  will be called the {\it mean bi-normal vector field} of $\alpha$.

The  eigenvalues $k_\alpha \leq  K_\alpha $ of the {\it Weingarten operator}  ${\mathcal W}_\alpha =-\Pi_{\alpha,t}DN_\alpha $ of ${\mathbb T}{\mathbb M}^2$ are called the  {\it principal mean curvatures} of $\alpha$.  The set of points ${\mathcal S}_{u,\alpha}$, where $k_\alpha =  K_\alpha$, will be called the {\it umbilic singularities} of $\alpha$.
 Outside   ${\mathcal S}_\alpha = {\mathcal S}_{n,\alpha}\cup{\mathcal S}_{u,\alpha}$, the {\it singular set} of $\alpha$,   are defined  the {\it  minimal}, $L_{m,\alpha}$,  and the {\it  maximal}, $L_{M,\alpha}$,   {\it principal mean line fields} of $\alpha$, which are  the eigenspaces of  ${\mathcal W}_\alpha$ associated respectively to  $k_\alpha$ and $ K_\alpha $.
The integral foliations of these line fields, which are of class $C^{r-2}$ on the complement of ${\mathcal S}_\alpha$, will be denoted by  ${\mathcal F}_{m,\alpha}$  and  ${\mathcal F}_{M,\alpha}$.

In a local chart $(u,v)$ the  principal mean curvatures lines  of $\alpha $ are   characterized as the solutions of the following  quadratic differential equation:
\begin{equation} \label{eq:lh}
(Fg_H-f_HG)dv^2+(Eg_H-e_HG)dudv+(Ef_H-Fe_H)du^2=0,\end{equation}
\noindent where $E=<\alpha_{u},\alpha_{u}>, \; F=<\alpha_{u},\alpha_{v}>, \; G=<\alpha_{v},\alpha_{v}>$ are the coefficients of the {\it first fundamental form} $I_\alpha =\alpha^{*}<\, , \,>$ and  $e_H=<\alpha_{uu},H_\alpha>,$ \;
$f_H=<\alpha_{uv},H_\alpha>,\; g_H=<\alpha_{vv},H_\alpha>$ the coefficients of the {\it second fundamental form relative to} $N_\alpha$, denoted $II_{N_\alpha}$, which have been   multiplied by $|H_\alpha|$ to remove the denominators.

The left hand member of equation (\ref{eq:lh}) is equivalent to the {\it Jacobian}, $J_\alpha$, of the quadratic  forms $I_\alpha$ and $II_{N_\alpha}$.

The {\it principal mean configuration} of $\alpha$ is defined by the quadruple
 ${\mathcal P}_\alpha =\{ {\mathcal S}_{n,\alpha} , {\mathcal S}_{u,\alpha}, {\mathcal F}_{m,\alpha},{\mathcal F}_{M,\alpha}\}.$

An immersion $\alpha$ in ${\mathcal I}^{r,s}$ is said to be {\it (r,s)-- principal mean curvature structurally stable} if it has a neighborhood $\mathcal V$ such that  for every $\beta$ in $\mathcal V$ there is a homeomorphism $h_\beta$ mapping ${\mathcal S}_{n,\beta}$ and ${\mathcal S}_{u,\beta}$ respectively  onto ${\mathcal S}_{n,\alpha}$ and  ${\mathcal S}_{u,\alpha}$ and
mapping the lines of the foliations ${\mathcal F}_{m,\beta}$ and  $ {\mathcal F}_{M,\beta}$ respectively  onto those of the foliations  ${\mathcal F}_{m,\alpha}$ and ${\mathcal F}_{M,\alpha}$. Denote by ${\mathcal E}^{r,s}$ the class of {\it (r,s)-- principal mean curvature structurally stable} immersions.

This global notion  can be localized at the {\it singularities}: $\{ {\mathcal S}_{n,\alpha} , {\mathcal S}_{u,\alpha}\}$ and at other
 invariant sets of the foliations, such as  some of the   {\it principal mean curvature cycles}, which are the periodic leaves  of the foliations.

A singularity of $\alpha$ is called {\it Darbouxian} if a) at all points of  the projective line over it, $dJ \neq 0$ and b)
the {\it Lie-Cartan Vector Field}, given locally by

\begin{equation}\label{eq:lc}
X_\alpha =(J_p,\,pJ_p,\,-(J_u+pJ_v)),
\end{equation}

\noindent with $ J =J_\alpha$ defined  in equation \ref{eq:lh},
is such that   along the projective line over the singularity has only hyperbolic equilibria  \cite{col}.  See Section \ref{sec:2} for more precise definition.

As in the standard   $\mathbb R ^3$ case, there are three Darbouxian types. In the case $D_1$ there is  only one hyperbolic saddle. In the case $D_2$ there are  three  hyperbolic singular points, one node and two saddles. In the case $D_3$ there are three hyperbolic saddle points. The subscript $i$ in $D_i$ denotes  the number of   {\it separatrices} reaching the singularity. Conditions on the third order jet of $\alpha$ at a singularity  to be Darbouxian  and which  discriminate its $D_i$ type will be given in Section \ref{sec:2}.   See Fig. \ref{fig:darboux}.

A {\it principal mean curvature cycle} $c$ %, that is a periodic integral curve
 of ${\mathcal F}_{m,\alpha}$  or  ${\mathcal F}_{M,\alpha}$,  is called {\it hyperbolic} if the derivative of its {\it first return} -- also called {\it holonomy} or {\it Poincar{\'e}} map,  $\pi_c$,  is different from $1$. An integral expression for  this derivative in terms of geometric curvature function along $c$ is given in Section \ref{sec:cyc}.

If the {\it limit set} of
a leaf of a principal mean foliation is contained
 in the set of  singular points and cycles it is said to  be  {\it limit set trivial}.

 The main result of this paper  can be stated now. The proof is outlined in Section \ref{sec:prf}.

\begin{theorem}\label{th:1}
Denote by  ${\Sigma}^{(r,s)}$ the class of immersions  which satisfy  the conditions on Darbouxian singularities, hyperbolic cycles, non-connection of singularity separatrices and triviality of the limit sets of all principal mean curvature lines. Then ${\Sigma}^{(r,s)}$ forms an open set in ${\mathcal I}^{r,s}$  and  it is contained in ${\mathcal E}^{r,s}$
for $r\geq 4$, $s\geq 3$.
\end{theorem}

The analysis of the density of  the class $\Sigma ^{(r,s)}$   in  ${\mathcal I}^{r,2}$ will be postponed to a forthcoming paper.
Meanwhile we will give some examples to confirm that $\Sigma ^{(r,s)}$  is not empty.

 \begin{example} {\bf a)} By composing with  the stereographic projection of ${\mathbb R}^3$ into ${\mathbb S}^3$ the  principal structurally immersions in ${\mathbb R}^3$  studied by   Gutierrez and Sotomayor  \cite{gs}, are obtained elements in $\Sigma ^ {(r,s)}$. See  \cite{lf}.

\noindent {\bf b)} Also, the subclass of Gutierrez and Sotomayor
of immersions in ${\mathbb R}^3$
 with non-vanishing mean curvature   is  contained inside $\Sigma ^{(r,s)}$.
\end{example}

Theorem \ref{th:1} partially extends a  result  of   Gutierrez and Sotomayor  \cite{gs}  for the    structural stability of principal configurations
on surfaces   in ${\mathbb R}^3$.  For other pertinent  extensions
to immersions of surfaces into ${\mathbb R}^4$,
see Section \ref{sec:5}.

\section{ Darbouxian Singularities }\label {sec:2}
 In a
Monge  chart $(u,v)$,  an immersion
$\alpha$ is expressed as follows:
\begin{equation}\label{eq:monge}\aligned
\alpha(u,v)=&(u,v,h_1(u,v),h_2(u,v)),\\
h_1(u,v)=&\frac{r_1}2 u^2+s_1uv+\frac{ t_1}2 v^2+\frac {a_1}6 u^3+\frac {d_1}2 u^2v+\frac {b_1}2 uv^2+\frac {c_1}6 v^3+O(4),\\
h_2(u,v)=&\frac{r_2}2 u^2+s_2 uv+\frac{ t_2}2 v^2+\frac {a_2}6 u^3+\frac {d_2}2 u^2v+\frac {b_2}2 uv^2+\frac {c_2}6 v^3+O(4).
\endaligned
\end{equation}
Let $B_1=(-\frac{\partial h_1}{\partial u},-\frac{\partial h_1}{\partial v}, 1,0)$, $B_2=\alpha_u\wedge \alpha_v\wedge B_1$. Write
$N_i=B_i/|B_i|, \; i=1,2$.
Clearly $<\alpha_u,N_i>=<\alpha_v,N_i>=<N_1,N_2>=0$   and

\begin{equation}\aligned  N_1(u,v)=&( -r_1u -s_1v+O(2), -s_1u -t_1v+O(2), 1, 0 ), \\
  N_2(u,v)=&( -r_2u -s_2v+O(2), -s_2u -t_2v+O(2), O(2), 1+O(2) )\endaligned
\end{equation}

The coefficients, $E,F,G$, of the first fundamental  form,  $I_\alpha$, induced by $\alpha$, and those   of  the second fundamental forms, $II_{i\alpha}$, relative to $N_i$, denoted  $e_i, f_i , g_i$, $i=1,\,2$, are calculated in the chart $(u,v)$ as follows:

\begin{equation}\label{eq:ff}
\aligned E(u,v)=& 1+O(2),\;\;\; F(u,v)=O(2),\;\;\; G(u,v)=1+O(2),\\
e_1(u,v)=& r_1+a_1u+d_1v+O(2),\;\;\; e_2(u,v)=  r_2+a_2u+d_2v+O(2), \\
f_1(u,v)=&s_1+d_1u+ b_1v+O(2), \;\; \; f_2(u,v)= s_2+d_2u+ b_2v+O(2)\\
g_1(u,v)=&t_1+b_1u+ c_1v+O(2), \;\;\; g_2(u,v)= t_2+b_2u+ c_2v+O(2)\endaligned
\end{equation}

The components
of $H= H_\alpha$ relative to a positive normal frame $N_1, N_2$ are given by
$H_1 = \frac{G\,e_1 -2F\,f_1+E\,g_1}{2(EG-F^2)},\, H_2=\frac{G\,e_2 -2F\,f_2 + E\,g_2}{2(EG-F^2)}$.  With the coefficients from equations   \ref{eq:ff}, in  the chart $(u,v)$, we obtain
\begin{equation}
\aligned H_1(u,v)=&\frac{t_1+r_1}2+\frac{a_1+b_1}2 u +\frac{c_1+d_1}2 v +O(2)\\
H_2(u,v)=&\frac{t_2+r_2}2+\frac{a_2+b_2}2 u +\frac{c_2+d_2}2 v +O(2).
\endaligned
\end{equation}
Thus
$H_\alpha =H_1N_1+H_2N_2$ can be written as follows:
\begin{equation}\label{eq:nh}
\aligned
H_{\alpha}=(& -( r_1 t_1+ r_1^2+ r_2 t_2 + r_2^2)\frac u2-( s_1 t_1+ s_1 r_1+ s_2 t_2+ r_2 s_2)\frac v2+O(2),\\
 -&( s_1 t_1+  s_1 r_1+ s_2 t_2 +r_2 s_2)\frac u2-(  t_1^2+ r_1 t_1+t_2^2+  r_2 t_2)\frac v2+O(2),\\
  & \frac{t_1+ r_1}2+(a_1+b_1)\frac u2+ ({c_1+d_1})\frac v2+O(2),\\
& \frac{t_2+ r_2}2+(a_2+b_2)\frac u2+ ({c_2+d_2})\frac v2+O(2))\endaligned
\end{equation}

Thus, in the chart $(u,v)$, a {\it normal singularity} located at $0$  is characterized by $(r_1+t_1 = 0,\; r_2+t_2=0)$.
Also the differential equation \ref{eq:lh}
of principal mean curvature lines  in the chart $(u,v)$ around   such a normal  singularity  is given by:
\begin{equation}\label{eq:sn}
 \aligned
   \{&- [s_1(a_1+b_1)+s_2(a_2+b_2)]\frac{u}{2}- [ s_1(c_1+d_1)+s_2(c_2+d_2)]
\frac{v}{2}+O_1\} dv^2 \\
 -& [(r_1(a_1+b_1)+r_2(a_2+b_2))u+  (r_1(c_1+d_1)+r_2(c_2+d_2)) v+ O_2]du dv + \\
   \{&[s_1(a_1+b_1)+s_2(a_2+b_2)]\frac{u}{2}+ [ s_1(c_1+d_1)+s_2(c_2+d_2)]\frac{v}{2}+O_3\}
  du^2=0 \\
   :=& -(\bar d u +\bar b v+O_1) dv^2  +[ \bar a u+  \bar c v+O_2]du dv+(\bar d u +\bar b v+O_3) du^2  =0
\endaligned
\end{equation}
\noindent where, $O_i=O(u^2+v^2)$.

 \begin{remark}\label{rem:rotn}
After an appropriate rotation in the frame $(u,v)$, it can be assumed that $\bar d=s_1(a_1+b_1)+s_2(a_2+b_2)=0$.
In fact, the equation \ref{eq:sn} in the  coordinates $(u_1,v_1)$, where  $u=\cos {\omega} u_1 +\sin {\omega} v_1, \;\; v=-\sin {\omega} u_1 +\cos {\omega} v_1$,
is given  by:
\begin{equation}\label{eq:sn1}
 \aligned
  -[{\bar d}_1& u_1 +{\bar b}_1 v_1] dv_1^2  +[ {\bar a}_1 u_1+  {\bar c}_1 v_1]du_1 dv_1
+[{\bar d}_1 u_1 +{\bar b}_1 v_1] du_1^2 + O(2)=0 \\
&\hskip -.8cm\text{ where,}\\
 {\bar  d}_1 =& \cos^3 {\omega} [\bar b \tan^3 {\omega} +(\bar c-\bar d)\tan^2{\omega}  -(\bar a+\bar b)\tan {\omega} +\bar d ]
\endaligned
\end{equation}

Solving the cubic equation $\bar d_1(\tan{\omega})=0$ the assertion  follows.
\end{remark}

A normal singularity is called {\it Darbouxian}  if

\noindent {\bf a)} $H_\alpha$ is transversal to the zero section of the normal bundle:

 $\bar a\bar b-\bar c \bar d = \bar a_1\bar b_1-\bar c_1 \bar d_1=
\frac 12 b_1(c_1+d_1)s_1r_1+\frac 12 b_2(c_2+d_2)s_2r_2+
 \frac 1 2[(a_1+b_1)(c_2+d_2) -a_2(c_1+d_1)]s_1r_2+
 \frac 12[(a_2+b_2)(c_1+d_1)-a_1(c_2+d_2)]s_2r_1 \ne 0$

and

\noindent {\bf b)} one of the following condition holds:
\begin{itemize}
\item[$D_1$)] $\bar d_1=0$,\;  $\bar c_1^2+4\bar b_1(\bar a_1+\bar b_1)<0,$

\item[$D_2$)] $\bar d_1=0$,\; $\bar c_1^2+4\bar b_1(\bar a_1+\bar b_1)>0, \; -1\ne \bar a_1/\bar b_1<0, $

\item[$D_3$)] $\bar d_1=0$,\;  $\bar a_1/\bar b_1>0. $
\end{itemize}

\begin{remark}
It can be shown that the conditions $D_i$  above are independent of the rotation performed to have $\bar d_1(\tan w)=0$ in Remark \ref{rem:rotn}.
\end{remark}

The differential equation \ref{eq:lh} of mean curvature lines  near the umbilic singularity $0$
     characterized by
 $   e_HG-g_HE= t_1^2+ t_2^2- r_1^2- r_2^2= 0\;\;\;$
and $\;\;\;\; f_HG-g_HF =s_1(t_1+r_1)+s_2(t_2+r_2)=0, $
 is given by:
\begin{equation}\label{eq:su}
\aligned
-&\{[  2r_1d_1  + s_1(a_2+b_2)]\frac u2+ [ 2 r_1b_1  + s_1(c_2+d_2)]\frac v2+O_1(2)\} dv^2 \\
+& \{[r_1( b_1 -  a_1)- r_2(a_2+b_2)] u+ [r_1(c_1-d_1)-r_2(c_2+d_2)] v+O_2(2)\} dudv \\
+&\{[  2r_1d_1  + s_1(a_2+b_2)]\frac u2+ [ 2 r_1b_1  + s_1(c_2+d_2)]\frac v2+O_3(2)\} du^2=0\\
 :=& -(\tilde d u + \tilde b v)dv^2 +[ \tilde a u+  \tilde c v]du dv+ (\tilde d u+  \tilde b v) du^2 + O(2)=0
\endaligned
 \end{equation}

\begin{remark}\label{rem:rotu}
As in the normal singularity case by an appropriate rotation in the plane
$(u,v)$ it can be assumed that $\tilde d= 2r_1 d_1  + s_1(a_2+b_2)=0$.
\end{remark}

A umbilic singularity    is called {\it Darbouxian}  if

\noindent {\bf a)}
${\mathcal W}_\alpha$,
 regarded  as a section
is transversal to the line bundle of diagonal operators; in terms of the coefficients defined in \ref{eq:su},  this transversality condition writes:
$ \tilde a\tilde b-\tilde c\tilde d=[b_1(b_1 -  a_1)+d_1(d_1 -  c_1)]r_1^2 +[d_1(c_2+d_2) -b_1(a_2+b_2)]r_1r_2+\frac 12[(a_2+b_2)(d_1 -  c_1)+(b_1 -  a_1)(c_2+d_2)]r_1s_2\ne 0$.
and

\noindent {\bf b)}  one of the following condition, expressed assuming the simplification in Remark \ref{rem:rotu},  holds:
\begin{itemize}
\item[$D_1$)] $\tilde d=0$,\;  $\tilde c^2+4\tilde b(\tilde a+\tilde b)<0,$
\item[$D_2$)] $\tilde d=0$, \; $\tilde c^2+4\tilde b(\tilde a+\tilde b)>0, \;; -1\ne \tilde a/\tilde b<0, $
\item[$D_3$)] $\tilde d=0$, \; $\tilde a/\tilde b>0. $
\end{itemize}

\begin{figure}[ht]
 \centerline{
  \epsfxsize=10cm\epsfysize=3cm
 \epsfbox{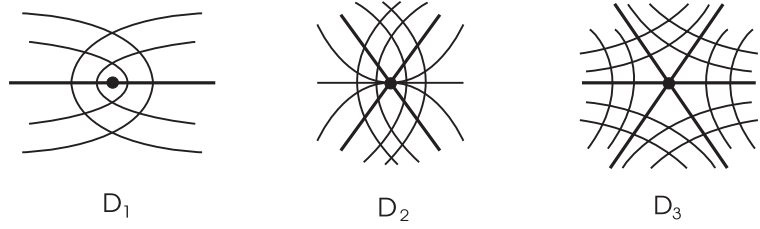 } }

\caption{Normal and Umbilic  Darbouxian Singularities and their Separatrices.
 \label{fig:darboux}  }
\end{figure}

The local behavior of the foliations ${\mathcal F}_{m,\alpha}$  and  ${\mathcal F}_{M,\alpha}$
near
singularities are as shown in the Fig. \ref{fig:darboux}.

 The conditions $D_i$ given above are similar to those obtained by  Gutierrez and Sotomayor  \cite{gs}   to characterize Darbouxian umbilic points of surfaces of $\mathbb R^3$. We have the following  correspondences with the $(a,\,b,\,c)$ notation of   \cite{gs} and   \cite{col}:  $\bar b_1=b, \;\bar c_1=c,\; \bar a_1=b -  a$, for normal singularities,  and  $\tilde b_1=b, \;\tilde c_1=c,\; \tilde a_1=b - a$, for umbilic singularities.

The proof of the local configurations for both cases is therefore the same as in   \cite{gs} and   \cite{col}, applied to the equilibria of the Lie-Cartan Vector Field \ref {eq:lc},  with $J_\alpha$ as in left hand members of the last equations in \ref{eq:sn} and   \ref{eq:su}.

\section{Hyperbolic Cycles}\label{sec:cyc}

Let $\alpha \in {\mathcal I}^{r,s} $ and suppose that $c$
is a regular  arc length parametrized curve in $\mathbb M ^2 \setminus {\mathcal S}_\alpha$.  Call  $t = c^\prime $ the
tangent vector field of $c$ and let  $T$ be  the unit vector field
along $c$  such that the tangent frame $\{ t,\ T \}$ is positive.

The equations of Darboux for the frame $\{t,\ T,\ N:=N_\alpha,\ B:=B_\alpha  \}$ along $c$ are given by:

\begin{eqnarray}\label{eq:d}
t^\prime=& k_gT + k N +  k_{B} B, \;\; T^\prime =& - k_g t -\tau_N N -  \tau B,
       \nonumber \\
N^\prime=&  - k t  + \tau_N T+ \tau_B B, \;\; B^\prime=&  - k_B t + \tau T - \tau_B N .
\end{eqnarray}

\begin{lemma}\label{lem:c}
Let $c$ be a minimal   principal mean  cycle of length $L$ of an immersed
surface $\mathbb M ^2$ in ${\mathbb R}^4$.

Then the expression

\begin{equation}\label{eq:ch}
\aligned \alpha(u,v) =&  c(u)  + vT(u) + [ K(u) \frac{v^2}2
+ a(u)\frac{v^3}6+ v^3 A(u,v)]N(u)\\
 +& [{\bar K}(u) \frac{v^2}2 +b(u)\frac{v^3}6+ v^3 {\bar A}(u,v)]B(u),\endaligned
\end{equation}

\noindent where $A(u,0)={\bar A}(u,0)=0$, defines  a L- periodic chart in a
neighborhood of $c$. \end{lemma}
\begin{proof} See     \cite{axial} and   \cite{gs}.
\end{proof}
With the notation in equations \ref{eq:d} and \ref{eq:ch}, follows that $c$ is a minimal principal mean  cycle if and only if the following holds  along it,
\begin{equation} \label{eq:k}
\tau_N\equiv 0,\; \; \; k_B+\bar K\equiv 0,\;\;\; K-k> 0.
\end{equation}

\begin{lemma}\label{lem:der}
Let $c$ be a  principal   mean   curvature cycle and consider a coordinate chart $(u,v)$
as in  Lemma \ref{lem:c}. Then the first derivative of the Poincar\'e
map $\pi$
of the principal cycle $c$ has the form

\begin{equation}\label{eq:pp}
ln \pi^\prime(0)= \int_0^L \frac{-[Ef_H -Fe_H]_v}{Eg_H -Ge_H} du,
\end{equation}
 \noindent where L is the length of the   principal
cycle, $E,\ F,\ G$ are the coefficients of the first fundamental
form and $e_H, \ f_H,\ g_H$ are the coefficients of the second
fundamental form with respect to the  normal vector field $
H_\alpha$ given in equation \ref {eq:lh}, calculated relative to
the chart $(u,v)$ defined by \ref{eq:ch}.
\end{lemma}

\begin{proof} The differential equation of principal mean  curvature lines
is given in equation \ref {eq:lh}.

 As  $F(u,0)= f_H(u,0)=0,$ the result follows by differentiating the
equation   above with respect to the initial condition $v_0$ --
thus getting the linear variational  equation--. Recall that
$\pi(v_0)=v(L,v_0)$ where $v(u,v_0)$ is the solution of this
equation with initial condition $v(0,v_0)=v_0$. The expression for
$\pi^\prime(0)$ in equation \ref{eq:pp}   follows from the
integration of the  linear variational  equation.
\end{proof}

The  calculation that follows culminates in an expression of  the
integral in equation \ref{eq:pp} in terms of the functions of the
arc length $u$ defined in equations \ref{eq:d},  \ref{eq:ch} and \ref{eq:k},
leading to the integral in Proposition \ref{pr:pp}.

\begin{equation} \label{eq:auv}
\aligned
\alpha_u (u,v) =& (1-k_g v-k_B A_2) t+(\tau_NA_1+\tau A_2) T \\
+&(\frac{\partial A_1}{\partial u}-\tau_N v-\tau_B A_2)N+( \frac{\partial A_2}{\partial u}-\tau  v-\tau_B A_1)B\\
\alpha_v(u,v)=& T+\frac{\partial A_1}{\partial v}N+\frac{\partial A_2}{\partial v}B\\
 A_1(u,v)=& K(u) \frac{v^2}2
+ a(u)\frac{v^3}6+ v^3 A(u,v) \\
A_2(u,v)=& {\bar K}(u) \frac{v^2}2 +b(u)\frac{v^3}6+ v^3 {\bar A}(u,v).
\endaligned
\end{equation}

Write $\alpha_u= x_1 t+x_2 T+x_3 N+x_4 B$ and $\alpha_v=T+y_1 N + y_2 B$.
Let $\bar N_1=(y_1 x_2-x_3)t -x_1y_1 T+x_1N +0B$ and
$\bar N_2=\alpha_u\wedge \alpha_v\wedge \bar N_1$.
Then it follows that $<\bar N_i, \alpha_u>=<\bar N_i,\alpha_v>=<\bar N_1,\bar N_2>=0$.
Direct calculations show that

\begin{equation}\label{eq:n2b}
\aligned \bar N_2=&(x_1x_2y_2-x_1x_4+x_1x_3y_1y_2-x_1 x_4y_1^2)t\\
+&
(-x_1^2y_2+x_2 x_3 y_1 y_2-x_2 x_4y_1^2-x_3^2y_2+x_3 x_4 y_1)T\\
 +&(-x_1^2 y_1y_2-x_2^2y_1y_2+x_2 x_4y_1+x_2x_3y_2-x_3x_4)N \\
+&(x_1^2+x_1^2y_1^2+x_2^2y_1^2-2x_2 x_3y_1+x_3^2)B.
\endaligned
\end{equation}

 Let $N_1(u,v)=\bar N_1/|\bar N_1| $ and $N_2(u,v)=\bar N_2/|\bar N_2|$  be orthonormal vector fields.

Straightforward calculations lead to:

\begin{equation}\label{eq:n1n2}
\aligned
N_1(u,v)=&[\tau_N v+O(2)] t+[ -K v +O(2)]T + [1+O(2)]N \\
N_2(u,v)=&[\tau  v+O(2)] t+[ -\bar K v +O(2)]T + [ O(2)]N +[1+O(2)]B
\endaligned
\end{equation}

From equation \ref{eq:auv} it follows that
\begin{equation} \label{eq:1f}
\aligned E(u,0)  = & G(u,0)=1,\; F(u,0)=0\\
E_v(u,0)=&-2k_g(u),\; F_v(u,0)=G_v(u,0)=0.\endaligned
\end{equation}
Also from equations \ref{eq:d}, \ref{eq:ch} and \ref{eq:n1n2}
 it follows that

\begin{equation} \label{eq:2f1} \aligned e_1(u,0)   =&k(u), \; \;f_1(u,0)=-\tau_N, \; \; g_1(u,0)=K(u)\\
(e_1)_v(u,0)=& -k_g(k+K)+\tau\tau_B -\tau_N^\prime\\
(f_1)_v(u,0)=& K^\prime-\bar K \tau_B -k_g\tau_N \\
(g_1)_v(u,0)=& a(u).\endaligned \end{equation}

Here $e_1(u,v)=<\alpha_{uu},N_1(u,v)>, \; f_1(u,v)=<\alpha_{uv}, N_1(u,v)>$ and $g_1(u,v)=<\alpha_{vv}, N_1(u,v)>.$

From   $e_2(u,v)=<\alpha_{uu},N_2(u,v)>, \; f_2(u,v)=<\alpha_{uv}, N_2(u,v)>$ and $g_2(u,v)=<\alpha_{vv},N_2(u,v)>$, it follows that:

\begin{equation} \label{eq:2f2}
\aligned e_2(u,0)  =&\,k_B(u), \;\; f_2(u,0)=-\tau(u), \;\; g_2(u,0)=\bar K(u)\\
(e_2)_v(u,0)=&\, -k_g(k_B+\bar K)  -\tau_N\tau_B-\tau^\prime \\
(f_2)_v(u,0)=& \,{\bar K}^\prime+  K \tau_B -k_g \tau \\
(g_2)_v(u,0)=& \, b(u).
\endaligned
\end{equation}

Define
\begin{equation} \label{eq:hi}
\tilde H_i=\frac{Eg_i-2f_iF+e_iG}{2(EG-F^2)},\;\;\; i=1, \; 2.\end{equation}

 Accordingly, the mean curvature vector writes ${ H}_\alpha(u,v)=\tilde H_1 N_1(u,v)+\tilde  H_2  N_2(u,v)$.

From equations \ref{eq:2f1},  \ref{eq:2f2} and \ref{eq:hi} it follows that:

\begin{equation}\label{eq:h12}
\aligned \tilde H_1(u,0)=&\frac{k+K}2, \;\;\;\;\tilde H_2(u,0)=\frac{k_B+\bar K}2=0\\
 (\tilde H_1)_v =& a(u)-k_g(K-k)+\tau\tau_B-\tau_N^\prime\\
2(\tilde H_2)_v =&b(u)-2k_g \bar K  -\tau_N\tau_B-\tau^\prime
\endaligned
\end{equation}

Therefore it follows from equations \ref{eq:n1n2} and \ref{eq:h12} that
 the  functions $e_H(u,v)=<\alpha_{uu},{  H}_\alpha(u,v)>, \; f_H(u,v)=<\alpha_{uv},
{ H}_\alpha(u,v)>$ and $g_H(u,v)=<\alpha_{vv},{ H}_\alpha(u,v)>$ evaluated at $v=0$ give:

 \begin{equation} \label{eq:efgH}
\aligned e_H(u,0)=& k\tilde H_1(u,0), \;\; f_H(u,0)=-\tau_N\tilde H_1=0, \;\; g_H(u,0)=K \tilde H_1(u,0)\\
(f_H)_v(u,0)=& \tilde H_1(u,0)(K^\prime +\tau_B \bar K )+(\tilde H_2)_v \tau.
\endaligned
\end{equation}

\begin{proposition}\label{pr:pp}
The first derivative of the Poincar\'e map of a minimal principal
cycle is given by

\begin{equation}\label{eq:dp}
 ln \pi^\prime(0)= - \int_0^L
          \frac{k^\prime}{K- k}du
         +\int_0^L  \frac{   k_B \;\tau_B   }
         {K- k}du - \int_0^L  \frac{ (\tilde H_2)_v (u,0)\; \tau }
         {\tilde H_1 (u,0)(K- k)}du
 \end{equation}
\end{proposition}
\begin{proof} It follows directly  from lemma \ref{lem:der},
 equations \ref{eq:h12} and \ref{eq:efgH} and integration by parts.
\end{proof}

\begin{remark}\label{rk:g}
In the last integral,  the expressions involving $\tilde H_1$ and $(\tilde H_2)_v$ can
 be further simplified using the equations in \ref{eq:h12}. Notice that this
 introduces $b(u)$ which however can itself be expressed in terms
of the three dimensional torsion of the  curve $v\to \alpha(u,v)$
in the $3-$space generated by $\{T(u), N(u), B(u)\}$.
\end{remark}

The next proposition shows how to deform an immersion making hyperbolic a cycle, under mild conditions.

\begin{proposition} \label{pr:def} Consider the
 one parameter family of immersions:

\begin{eqnarray}\label{eq:defor} \alpha_{\epsilon}(u,v) &=  \alpha(u,v) + \epsilon \delta(u)m(v)\frac{v^3}6 B(u)
\end{eqnarray}

\noindent where $m(v)=1$ in neighborhood of $v=0$, with small support and $\delta > 0$.

If      $\tau\ne 0$,
then $c$ is a  hyperbolic principal cycle for all immersions
$\alpha_{\epsilon},\; \epsilon>0 $ small.
\end{proposition}
\begin{proof} Along $c$   the deformation $\alpha_\epsilon$ given by equation \ref{eq:defor} has the same
 second order jet as that  of $\alpha$.    It follows that $c$ is also an arc length parametrized  minimal principal mean  curvature cycle for $\alpha_\epsilon$. In the integral  expression \ref{eq:dp}  for the   derivative of the Poincar\'e map it follows that $(\tilde H_2)_v(u,0,\epsilon)= (\tilde H_2)_v(u,0)+\epsilon \delta(u)$ while  all the other functions involved are independent on $\epsilon$. Therefore, after a direct calculation, it follows that
$\frac{\partial}{\partial \epsilon}(ln \pi_{\epsilon} ')(0)\mid_{\epsilon =0}
       = -\int_0^L   \frac{\tau \delta}
         {\tilde H_1(K- k)}du$, which is positive taking  $ \delta =  -\tau\tilde H_1$.
 \end{proof}

\section {Outline of the Proof of Theorem \ref {th:1} } \label{sec:prf}

Once the hypotheses on $\alpha$ are expressed in  the
Projective Tangent bundle of $\mathbb M ^2$ and identified
with those for the quadratic equation \ref{eq:lh} which,
in turn, amount to the hyperbolicity of equilibria  and
 periodic orbits of the   Lie-Cartan Line Field,  locally
 expressed by $X_\alpha$ in \ref{eq:lc}, the similarity with
the case of    principal line fields  dealt with in   \cite {gs}
and   \cite {col} becomes evident. In fact,  the construction and
continuation   to a small neighborhood ${\mathcal V}(\alpha)$
of $\alpha$ of the canonical regions follow also from the openness
 and unique continuation, for $\beta$ near $\alpha$, of the
singularities (and their separatrices and parabolic sectors)
and   of cycles (and their local invariant manifolds), due to
the hyperbolicity of these elements in the field $X_\alpha$. This
leads to the openness of $\Sigma^{r,s}$ and gives the uniqueness
of the correspondence between singularities, normal and umbilic,
 separatrices, cycles for both minimal and maximal foliations involved
 and their intersections for ${\mathcal P}_\alpha$ and
 ${\mathcal P}_\beta$. The extension of this correspondence to
 define a topological equivalence homeomorphism $h_\beta$, is
 carried out as in the case  of principal configurations   \cite{gs}.

\section{Concluding Remarks and Related Problems}\label {sec:5}

The study   of the bending of a surface immersed in
${\mathbb R} ^4$, focusing  the stability  properties of the integral foliations
  defined by  geometric properties related to certain normal line
fields, has a rich background. The approach  and pertinent results
 presented here should be considered in the perspective of previous
achievements. A concise discussion follows.

The  {\it axial configuration} of Garcia and Sotomayor  \cite{axial}, for which
the normal line fields are those of the {\it principal axes} of the
{\it ellipse of curvature}, may be richest of all them.
In fact, it  leads to fields of {\it tangent crosses} rather than to
 tangent line fields. This theory, when restricted to a surface immersed in
${\mathbb R} ^3$, reduces to both  the standard {\it principal}
and that of the {\it arithmetic mean}  \cite{amean} configurations.

Garcia and Sotomayor  \cite{gas}  have studied principal cycles  of immersions
of surfaces in a three dimensional Riemannian manifold. The expression of the derivative of the return map should be compared  with that of equation \ref{eq:dp}.

 By taking the normal line field to be an arbitrary unit vector field, $\nu$,
 Garcia and S\'anchez  \cite{gsb} have obtained   an integral
expression for the first derivative of the return map associated to
 a principal cycle.

 Mello  \cite{lf}  has considered the
tangent line fields defined by the property of having their
  normal curvature  vector parallel to $H_\alpha$. The approach  of the present paper is in between this and  the previous one.

The consideration of other geometric normal vector fields such as $\nu = B_\alpha$, the {\it bi-normal},  instead
 of $H_\alpha$  in the present paper, may be also of  interest.

The $C^2$ density of the  { \it limit set triviality}  condition  seems
to be most difficult problem left open here;
see Theorem \ref{th:1}. This problem is also present  and, as far as
we know, still open for the previous approaches mentioned above.

Other direction  of  research, though not directly related to stability,
 emerges  with the evaluation of the {\it Index} of an isolated
singularity of ${\mathcal P}_\alpha$.   This is related to the upper
bound  $1$ for the umbilic index on surfaces in  ${\mathbb R} ^3$,
connected  to   deep problems around the   {\it Carath\'eodory
Conjecture}.  See Smyth and Xavier  \cite{sx} and Ivanov  \cite{iva}.

Gutierrez and S\'anchez  \cite{gusa} have shown that this bound does not hold for the $\nu$ approach. The
case of $\nu=H_\alpha$ presented here  contrasts with the flexibility  in the case of arbitrary $\nu$.  The question of the upper bound of  the index of a singularity seems  more difficult to   analyze in the present case.

\vspace{1cm}

{\small
 \noindent \address{Instituto de Matem\'{a}tica e
Estat\'{\i}stica \\
Universidade Federal de Goi\'as\\
 Caixa Postal 131\\
 74001-970  Goi\^ania, GO, Brasil\\
 E-mail: ragarcia@mat.ufg.br }

\vspace{1cm}

 \noindent \address{Instituto de Ci\^encias\\
 Universidade Federal de Itajub\'a\\
 37500-903 Itajub\'a, MG, Brasil\\
 E-mail: lfmelo@unifei.edu.br}

\vspace{1cm}
 \noindent\address{Instituto de Matem\'{a}tica e
Estat\'{\i}stica\\
  Universidade de S\~{a}o Paulo\\
    Rua do Mat\~{a}o 1010,
Cidade Universit\'{a}ria\\
    05508-090  S\~{a}o Paulo, SP, Brasil\\
  E-mail: sotp@ime.usp.br}
}

\end{document}